\documentclass{article}
\usepackage[utf8]{inputenc}
\usepackage{amsthm, amsfonts, amssymb, latexsym, amsmath}
\usepackage{xcolor}
\usepackage{ulem}

\theoremstyle{definition}

\usepackage{refcheck}
\begin{document}
\title{\bf  Some interesting
 number theory problems}
\author{  Wenpeng Zhang$^{1, 2}$\\
{\small 1.\    School of Data Science and Engineering}\\
{\small  Xi'an Innovation College of Yan'an University, Xi'an, Shaanxi, P. R. China}\\
{\small 2.\ School of Mathematics, Northwest University, Xi'an, Shaanxi, P. R. China}\\
{\small  Correspondence: wpzhang@nwu.edu.cn}}
\date{}
\maketitle
\baselineskip 16pt \begin{center}
\begin{minipage}{120mm}
{\bf Abstract:} {\small The main purpose of this paper is to propose some interesting number theory problems related to the Legendre's symbol and the two-term exponential sums.} 
\end{minipage}
\end{center}

\section{ Several problems related to Legendre's symbol. }
For any positive integer $q>1$, the classical Kloosterman's sums $S(m, n; q)$ is defined as follows:
$$
 S(m, n; q) = \mathop{{\sum}'}_{a=1}^{q} e\left(\frac{ma +n\overline{a}}{q}\right),
 $$
where $m$ and $n$ are any integers, $\displaystyle \mathop{{\sum}'}_{a=1}^{q}$ denotes the summation over all $1\leq a\leq q$ such that $(a, q)=1$, $\overline{a}$ denotes $a\cdot \overline{a}\equiv 1 \bmod q$, $e(y)=e^{2\pi iy}$ and $i^2=-1$.

This sum was introduced by H. D. Kloosterman [1] to study the square partitions of an integer. Thereafter, it plays a very important role in the research of analytic number theory, many important number theory problems are closely related to it.  Because of this, some scholars have studied the various elementary properties of $ S(m, n; q)$, and obtained a series of important results. For example,  H. Sali$\acute{e}$ [3] (or see H. Iwaniec [2]) proved that for any odd prime $p$, one has the identity
 $$
 \sum_{m=0}^{p-1}\left|\sum_{a=1}^{p-1}e\left(\frac{a+m\overline{a}}{p}\right)\right|^4= 2p^3-3p^2-3p.
 $$

 W. P. Zhang [4] has used the elementary method to obtain a generalized result. That is, for any integer $n$ with $(n, q) = 1$, one has the identity
$$
\sum_{m=1}^q\left|\mathop{{\sum}'}_{a=1}^{q} e\left(\frac{ma+n\overline{a}}{q}\right)\right|^4=3^{\omega(q)}q^2\phi(q)\prod_{p\| q}\left(\frac{2}{3}-\frac{1}{3p}-\frac{4}{3p(p-1)}\right),
$$
where $\phi(q)$ is Euler function, $\omega(q)$ denotes the number of all different prime divisors of q, $p\|q$ denotes the product over all prime divisors of $q$ with $p \mid q$ and $p^2\nmid q$.

In all results on Kloosterman's sums, perhaps the most property is the upper bound estimate for $ S(m, n; q)$ (see S. Chowla [5] or T. Estermann [6]). That is,
 $$
 \mathop{{\sum}'}_{a=1}^{q} e\left(\frac{ma+n\overline{a}}{q}\right)\ll \left(m, n, q\right)^{\frac{1}{2}}\cdot d(q)\cdot q^{\frac{1}{2}},
 $$
 where $d(q)$ is the Dirichlet divisor function, and  $(m, n, q)$ denotes the greatest common factor of $m,\ n$ and $q$.

There are so many results for $S(m, n; q)$ that we would not list them all here, but some of them can be found in references [7]--[11].

However, it is worth mentioning that W. P. Zhang, L. Wang and X. G. Liu [12] studied the calculating problem of the fourth power mean
 \begin{eqnarray*}
 \sum_{m=0}^{p-1}  \left| \sum_{a=1}^{p-1}   e\left(\frac{ma^k+\overline{a}}{p}\right)\right|^4,
 \end{eqnarray*}
 and proved the following conclusion:
\begin{eqnarray*}
&&\sum_{m=0}^{p-1}  \left| \sum_{a=1}^{p-1}   e\left(\frac{ma^2+\overline{a}}{p}\right)\right|^4\nonumber \\
&=& \left\{
\begin{array}{ll}  \displaystyle 2p^3-6p^2-5p+2\left(\frac{3}{p}\right)\cdot p^2- p^2\cdot \sum_{c=1}^{p-1}\left(\frac{c+1+\overline{c}}{p}\right) &\textrm{ if $p\equiv 3\bmod 4$;}\\ \displaystyle 2p^3-10p^2-9p-2\left(\frac{3}{p}\right)\cdot p^2+ p^2\cdot \sum_{c=1}^{p-1}\left(\frac{c+1+\overline{c}}{p}\right) &\textrm{ if $p\equiv 1\bmod 4$,}
\end{array}\right.
\end{eqnarray*}
where $p>3$ is a prime, and $\left(\frac{*}{p}\right)$ denotes the Legendre's symbol modulo $p$.

Recently, S. S. Ning and X. X. Wang [13] also studied the same problem, and obtained  a different conclusion. That is, they proved that

\begin{eqnarray*}
&&\sum_{m=0}^{p-1}  \left| \sum_{a=1}^{p-1}   e\left(\frac{ma^2+\overline{a}}{p}\right)\right|^4\\
&=& \left\{
\begin{array}{ll}  \displaystyle 2p^3-4p^2+2\cdot p^2\left(\frac{3}{p}\right)-5p+ p^2\cdot C(p) &\textrm{ if $p\equiv 3\bmod 4$;}\\ \displaystyle 2p^3-8p^2-2\cdot p^2\left(\frac{3}{p}\right)-9p+ p^2\cdot C(p) &\textrm{ if $p\equiv 1\bmod 4$,}
\end{array}\right.
\end{eqnarray*}
where $C(p)= \displaystyle \sum_{b=1}^{p-1} \left(\frac{\left(b^2+1\right)\left(b^2+4b+1\right)}{p}\right)$.

From these two different conclusions we may immediately deduce the following:

{\bf Corollary.} For any odd prime $p$, we have the identity
\begin{eqnarray}
\left(\frac{-1}{p}\right)\sum_{c=1}^{p-1}\left(\frac{c^3+c^2+c}{p}\right)-\sum_{b=1}^{p-1}\left(\frac{\left(b^2+1\right)\left(b^2+4b+1\right)}{p}\right)=2.
\end{eqnarray}

That is to say, the difference between two character sums of different polynomials is a constant $2$. Of course, this corollary is obtained by using the fourth power mean of the generalized Kloostermann sums as a bridge. By numerical testing (for all primes $3\leq p< 200$), the identity (1) is correct. Our problems are as follows:
\bigskip

(A). \ Give a direct elementary proof for formula (1).
\bigskip

(B).\ Whether there are infinitely many pairs of fundamentally different integer coefficients polynomials $f(x)$ and $g(x)$ (That is, $\left(\frac{f(x)}{p}\right)\neq\left(\frac{g(x)}{p}\right)$) such that
\begin{eqnarray}
\sum_{x=1}^{p-1}\left(\frac{f(x)}{p}\right)-\sum_{x=1}^{p-1}\left(\frac{g(x)}{p}\right)=c,
\end{eqnarray}
where $c$ is a fixed constant.
\bigskip

(C).\ What are the common properties between the polynomials $f(x)$ and $g(x)$ that satisfy identity (2).

\bigskip

(D). \ Whether the values of $c$ can only be $0$ or $2$?

\section{A problem related to exponential sums. }

Let $q\geq3$ be an integer. For any integers $m$ and $n$, the two-term exponential sums $S(m, n, k; q)$ is defined as follows:
\begin{eqnarray}
S(m, n, k; q)=\sum_{a=1}^{q}e\left(\frac{ma^k+na}{q}\right),
\end{eqnarray}
where $e(y)=e^{2\pi iy}$ and $i^2=-1$.

If $k=2$ and $n=0$ in (3), then the two-term exponential sums becomes the quadratic Gauss sums. These sums are crucial to analytic number theory and have a tight relationship to well-known number theory problems like the Waring problem and Goldbach conjecture. Therefore, many academics have investigated the various properties of the two-term exponential sums, and have developed several significant results. In general, people pay attention to the research content of the two-term exponential sums in two aspects: One is its upper bound estimate, and the other is its power mean. In fact, if $q=p$ is a prime, then for any fixed positive integer $k$, from A. Weil's classical work [14] one can deduce the best upper bound estimate
\begin{eqnarray*}
\sum_{a=1}^{p}e\left(\frac{ma^k+na}{p}\right)\ll p^{\frac{1}{2}}.
\end{eqnarray*}

As for the power mean calculating problem of $S(m, n, k; q)$, there are many interesting results.  For example,  H. Zhang and W. P. Zhang [15] proved that for any odd prime $p$, one has
\begin{eqnarray*}
\sum_{m=1}^{p-1}\left|\sum_{a=0}^{p-1}\
e\left(\frac{ma^3+na}{p}\right)\right|^4 =
 \left\{
\begin{array}{ll} 2p^3-p^2 &\textrm{ if\  $3\nmid p-1$;} \\   2p^3-7p^2 &\textrm{ if\  $3|
p-1$,}\end{array}\right.
\end{eqnarray*}
where $n$ represents any integer with $(n, p)=1$.

W. P. Zhang and D. Han [16] obtained the identity
$$
\sum_{a=1}^{p-1}\left|\sum_{n=0}^{p-1}
e\left(\frac{n^3+an}{p}\right)\right|^6= 5p^4-8p^3-p^2,
$$
where $p$ denotes an odd prime with $3\nmid (p-1)$.

W. P. Zhang and Y. Y. Meng [17] proved that
\begin{eqnarray*}
\sum_{m=1}^{p-1}\left|\sum_{a=0}^{p-1}e\left(\frac{ma^3+na}{p}\right)\right|^6=\left\{
\begin{array}{ll}\displaystyle 5p^3(p-1) &\textrm{ if $p\equiv 5\bmod 6$;} \\  \displaystyle 5p^4-23p^3 -d^2 p^2 &\textrm{ if $p\equiv1\bmod6$,}
\end{array}\right.
\end{eqnarray*}
where $4p=d^2+27b^2$, and $d$ is uniquely determined by $d\equiv 1\bmod 3$.

T. T. Wang and W. P. Zhang [18] obtained the eighth power mean of $S(m, 1, 3; p)$ and proved the identity
\begin{eqnarray*}
\sum_{m=1}^{p-1}\left|\sum_{a=0}^{p-1}e\left(\frac{ma^3+na}{p}\right)\right|^8=\left\{
\begin{array}{ll}\displaystyle 7(2p^5-3p^4) &\textrm{ if\quad $p\equiv5\bmod6$;} \nonumber\\  \displaystyle 14p^5-75p^4-8p^3d^2 &\textrm{ if\quad $p\equiv1\bmod6$.}\nonumber
\end{array}\right.
\end{eqnarray*}

Recently, my PhD student Wang Li proposed the following:

{\bf Conjecture.}  Let $p$ be an odd prime. Then for any positive integer $k$, we have the asymptotic formula
\begin{eqnarray*}
\sum_{m=0}^{p-1}\left|\sum_{a=0}^{p-1}e\left(\frac{ma^3+a}{p}\right)\right|^{2k}= \frac{1}{k+1}\binom{2k}{k}\cdot p^{k+1} + O\left(p^{k+\frac{1}{2}}\right).
\end{eqnarray*}

 It is not difficult to see from [16]-[18] that this conjecture is correct when $k=1, \ 2,\ 3$ and $4$. In an unpublished paper, the author also proved that for $k=5$ and $6$, the conjecture  is also true.

\end{document}